\documentclass[twoside]{siamltex}

\usepackage[final]{graphicx}
\usepackage{amssymb}
\usepackage{amsmath}

\newcommand{\weg}[1]{ }
\newcommand{\Tt}[1]{\mathbf{#1}}
\def \dt    { h }

\def \Nxi     {N_\xi}

\title{Transparent boundary conditions based on the Pole Condition for
time-dependent, two-dimensional problems}

\author{Daniel Ruprecht\footnotemark[1] \and Achim Sch\"adle\footnotemark[2] \and Frank Schmidt\footnotemark[3]}

\begin{document}

\maketitle

\renewcommand{\thefootnote}{\fnsymbol{footnote}}
\footnotetext[1]{Institute of Computational Science, USI Lugano, CH-6904 Lugano, Switzerland.\\~E-mail: 
{\tt daniel.ruprecht@usi.ch}. Supported by the Swiss HP2C initiative.}
\footnotetext[2]{Mathematisches Institut, Heinrich-Heine-Universit\"at, D-40255 D\"usseldorf, Germany.\\~E-mail:
{\tt schaedle@am.uni-duesseldorf.de}}
\footnotetext[3]{ZIB Berlin, D-14195 Berlin, Germany. \\~E-mail: 
{\tt frank.schmidt@zib.de}. Supported by the DFG Research Center \textsc{Matheon} 
"Mathematics for key technologies" in Berlin.}

\renewcommand{\thefootnote}{\arabic{footnote}}

\begin{abstract}
The pole condition approach for deriving transparent boundary
conditions is extended to the time-dependent, two-dimensional case.
Non-physical modes of the solution are identified by the position of
poles of the solution's spatial Laplace transform in the complex
plane. By requiring the Laplace transform to be analytic on some
problem dependent complex half-plane, these modes can be
suppressed. The resulting algorithm computes a finite number of
coefficients of a series expansion of the Laplace transform, thereby
providing an approximation to the exact boundary condition. The
resulting error decays super-algebraically with the number of
coefficients, so relatively few additional degrees of freedom are
sufficient to reduce the error to the level of the discretization
error in the interior of the computational domain. The approach shows
good results for the Schr\"odinger and the drift-diffusion equation
but, in contrast to the one-dimensional case, exhibits instabilities
for the wave and Klein-Gordon equation. Numerical examples are shown
that demonstrate the good performance in the former and the
instabilities in the latter case.
\end{abstract}


\section{Introduction}\label{Sec:Intro}
Transparent boundary conditions (TBCs) are required whenever a problem
is posed on a domain that has to be truncated in order to become
numerically treatable, either because it is unbounded or too large to
compute solutions in a reasonable amount of time. Usually, TBCs have
to avoid reflections at the artificial boundary, although more complex
situations can arise, for example if inhomogeneities are present in
the truncated part. Exact TBCs are typically non-local in time and
space and suitable approximations have to be derived in order to be
able to efficiently compute numerical solutions to the truncated
problem. The study of this type of boundary conditions started in the 
1970s, see the paper of E. L. Lindman~\cite{Lindman75} and references
given there. In their seminal paper~\cite{EngquistMajda} Engquist and Majda 
devised a general strategy for the derivation of approximate
TBCs. Comprehensive overviews of the subject can be found, for
example, in~\cite{AntoineEtAl2008,Givoli91,Givoli2004,Hagstrom99,Tsynkov98}.

The pole condition approach for the derivation of TBCs was introduced
in a first version in~\cite{Schmidt1995, Schmidt1997} for
time-dependent Schr\"odinger-type equations, later in~\cite{Hohage03a,
Hohage03b,Schmidt02H} for time-harmonic scattering problems. It was
further explored
in~\cite{GanderSchaedle2010,HohageNannen2009,Nannen2008,SchmidtEtAl2008}.
An alternative formulation of the pole condition is presented
in~\cite{NannenSchaedle:11}, which provides a noticeably simplified
implementation and is also used in the present paper. A comparison of
different techniques to derive TBCs for Schr\"odinger's equation can
be found in~\cite{AntoineEtAl2008}, finding the pole condition to be
one of the most efficient.  In~\cite{RuprechtEtAl2008}, the pole
condition approach is adopted for a larger class of time-dependent
problems, showing good performance for different types of partial
differential equations (PDEs) ranging from Schr\"o\-ding\-er's
equation, the heat equation to wave and Klein-Gordon
equation. However, the experiments involved only one-dimensional or
two-dimensional wave-guide geometries.  The present paper extends this
approach to the fully two-dimensional case and investigates its
performance through numerical experiments. While the very good
performance of the pole condition is confirmed in the two-dimensional
case for Schr\"odinger's equation and the drift-diffusion equation,
instabilities are found for the wave equation.

As the infinite element method, see~\cite{Astley2000}, the pole
condition does not truncate the exterior domain at some finite
length. Nevertheless, the finite number of expansion coefficients of
the Laplace transform also results in some form of truncation and the
pole condition realizes a radiation boundary condition at the boundary
of the interior domain and does not aim at providing a meaningful
solution in the exterior. In some special cases, see
\cite{Zschiedrich03}, the pole condition is closely related to the
perfectly matched layer approach introduced in~\cite{Berenger94}, but
as it does not require complex coordinate stretching, the pole
condition provides a more general framework. Note that in contrast to
other approaches to TBC involving Laplace transforms, for example
\cite{Alpertetal98b}, the pole condition applies the Laplace transform
in space and not in time.

The class of problems considered are, as in~\cite{RuprechtEtAl2008},
initial value problems for linear PDEs of the form
\begin{equation}
  p(\partial_{t}) u(t,\Tt{x}) = 
  c^{2} \Delta u(t,\Tt{x}) - \Tt{d} \cdot \nabla u(t,\Tt{x})
  - k^2 u(t,\Tt{x})
  \ \text{for} \ \Tt{x} \in\mathbb{R}^{2}, \ t \ge 0.
  \label{eq:u}
\end{equation}
Included here are the Klein-Gordon equation for $p(\partial_{t})
= \partial_{tt}$ and $\Tt{d} = (0,0)^{\rm T}$, the drift-diffusion
equation for $p(\partial_{t}) = \partial_{t}$ and $k=0$, the heat
equation for $p(\partial_{t}) = \partial_{t}$ and $\Tt{d} = (0,0)^T$, $k=0$
and finally Schr\"odinger's equation for $p(\partial_{t}) =
i \partial_{t}$ and $\Tt{d}=(0,0)^{\rm T}$, $k=0$.
Equation~\eqref{eq:u} is to be solved on a finite computational domain
$\Omega \subset \mathbb{R}^{2}$ with some boundary condition
$\Tt{B}(u) = 0 \ \text{on} \ \partial \Omega$, such that on the domain
$\Omega$ the solution of the initial boundary value problem
approximates the solution of the unrestricted initial value problem.

If the support of the initial value $u(0, \Tt{x})$ is a subset of
$\Omega$ and the exterior domain is homogeneous, in the linear case
the boundary condition has to suppress all modes traveling from the
exterior $\mathbb{R}^{2} \backslash \Omega$ into the computational
domain. Section~\ref{sec:poleCond} illustrates the main concept of
the pole condition by means of a simple one-dimensional
example. Section~\ref{sec:Disc} introduces the details of the
discretization employed in the two-dimensional case and
Section~\ref{sec:numExam} shows several numerical examples.

\section{Pole condition}\label{sec:poleCond}
This section provides a brief sketch of the key idea of the pole
condition. Denote the Laplace transform of some function $f$
along some (spatial) coordinate $r$ by
\begin{equation}
  \mathcal{L}(f)(s) = \int_{0}^{\infty} \exp(-s r) f(r) \ dr.
\end{equation}
The pole condition exploits the identity
\begin{equation}
\label{eq:lapTransExp}
  \exp(a r) \stackrel{\mathcal{L}}{\mapsto} \frac{1}{s - a},
\end{equation}
that is a mode with phase $a$ in physical space corresponds to a pole
of the Laplace transform located at $a$. The poles of the
Laplace transform of the solution are decomposed into poles
corresponding to incoming and outgoing modes or, more generally, into
poles corresponding to physical and non-physical modes. If the
locations of the poles in the complex plane corresponding to these
modes can be separated by a line, one can decompose the complex plane
into a half-plane $\mathbb{C}_{\rm in}$ containing all incoming modes
and a half-plane $\mathbb{C}_{\rm out}$ containing all outgoing
modes. Note that these half-planes depend on the equation at hand:
Table~\ref{tab:Cin} quotes the regions corresponding to the equations
mentioned above as derived in~\cite{RuprechtEtAl2008}. For a given
$\mathbb{C}_{\rm in}$, the pole condition is then defined as follows:
\begin{definition}\label{def:poleCond} 
  Let $u(t,r)$ be a function depending on time $t$ and some spatial
  coordinate $r$. Denote its Fourier transform in $t$ by $\hat{u}$ and
  the dual variable to $t$ by $\omega$. Then $u$ satisfies the pole
  condition, if $U(\omega, s) := \mathcal{L}(\hat{u}(\omega,\cdot))(s)$ 
  has an analytic extension to $\mathbb{C}_{\rm in}$ for every $\omega$.
\end{definition}

\begin{table}
  \caption{Region $\mathbb{C}_{\rm in}$ for different equations as derived 
    in~\cite{RuprechtEtAl2008}.\label{tab:Cin}}
  \centering
  \begin{tabular}{| c | c | c |}
    \hline
    Equation & Parameters in \eqref{eq:u} & $\mathbb{C}_{\rm in}$ \\ \hline
    Schr\"odinger equation & $p(\partial_{t}) = i \partial_{t}$, $\Tt{d}=0, k = 0$ & 
    $\left\{ z \in \mathbb{C} : \textrm{Re}(z) > -\textrm{Im}(z) \right\}$ \\
    Drift-diffusion equation & $p(\partial_{t}) = \partial_{t}$, $k=0$ 
    & $\left\{ z \in \mathbb{C} : \textrm{Re}(z) > 0  \right\}$ \\
    Wave equation & $p(\partial_{t}) = \partial_{tt}$, $\Tt{d}=0, k=0$ 
    & $\left\{ z \in \mathbb{C} : \textrm{Im} < 0 \right\}$ \\
    Klein-Gordon equation  & $p(\partial_{t}) = \partial_{tt}$, $\Tt{d}=0$ 
    & $\left\{ z \in \mathbb{C} :  \textrm{Im} < 0 \right\}$ \\ \hline
  \end{tabular}
\end{table}
To illustrate this concept, consider the one-dimensional wave equation
on a semi-infinite interval
\begin{equation}
  \label{eq:examWE}
  \partial_{tt} u(x,t) = \partial_{xx}u(x,t), \quad x \in \Omega = [-1, \infty)
\end{equation}
and  assume that a boundary condition at $x=0$ is sought such that the
solution of \eqref{eq:examWE} coincides with the solution on the restricted
domain $[-1,0]$. Here, the $r$ from Definition \ref{def:poleCond} is identical to the spatial coordinate $x$. Inserting an ansatz 
\begin{equation}
  u(x,t) = \exp( - i \omega t ) \exp(i k x)
\end{equation}
into~\eqref{eq:examWE} yields the dispersion relation
\begin{equation}
  \omega = \pm k
\end{equation}
and assuming $\omega > 0$ without loss of generality yields solutions of the form
\begin{equation}
  \label{eq:WEsol}
  u(x,t) = c_{1} \exp\left( - i \omega t \right) \exp\left( i \omega x
\right) + c_{2} \exp\left( - i \omega t \right) \exp\left(- i \omega x
\right),
\end{equation}
where the first term corresponds to the positive branch of the
dispersion relation and is rightward moving while the second term
corresponds to the negative branch and is leftward
moving. Let the non-physical modes in this example be the modes traveling
leftwards from $(0,\infty)$ into the interval $[-1, 0]$. The pole condition then has 
to suppress the pole corresponding to the second term in~\eqref{eq:WEsol}.

In order to point out the connection between the two modes
in~\eqref{eq:WEsol} and their corresponding poles, we derive an
equation for $U = \mathcal{L}(u)$ from \eqref{eq:examWE}. The Laplace transform
satisfies the identity
\begin{equation}
	\label{eq:lapTransId}
  \mathcal{L}(\partial_{xx}f)(s) = s^{2} \mathcal{L}(f)(s) - s f_{0} - f_{0}',
\end{equation}
where $f_{0}$ and $f_{0}'$ denote the Dirichlet and Neumann data at
$x=0$. Further, as in Definition~\ref{def:poleCond}, denote by $U$ the
function obtained by applying to $u$ Fourier transform in time and
Laplace transform in space. Using~\eqref{eq:lapTransId}, we obtain
from~\eqref{eq:examWE} the equation
\begin{align}
  -\omega^{2} U(\omega, s) &= s^{2} U(\omega, s) - s \hat{u}_{0}(\omega) -
  \hat{u}_{0}'(\omega) \nonumber \\ 
  \Rightarrow \quad U(\omega, s) &=
  \frac{1}{2} \frac{\hat{u}_{0} - (i/\omega) \hat{u}_{0}'}{s - i \omega} +
  \frac{1}{2} \frac{\hat{u}_{0} + (i/\omega) \hat{u}_{0}'}{s + i \omega},
  \label{eq:Udecomp}
\end{align}
where $\hat{u}_{0}(\omega)$, $\hat{u}_{0}'(\omega)$ are the Fourier transforms of
the time-dependent Dirichlet and Neumann data $u_{0}(t)$, $u_{0}'(t)$
at $x=0$. By~\eqref{eq:lapTransExp}, the first term with pole at $i\omega$ 
corresponds to the physically correct rightward propagating
mode with coefficient $c_{1}$ in~\eqref{eq:WEsol}, the second term
with pole at $ -i \omega$ to the non-physical leftward traveling
mode with coefficient $c_{2}$. In order to exclude the non-physical
mode, one could for example set $\mathbb{C}_{\rm in} = \left\{ z \in
\mathbb{C} : \textrm{Im}(z) < 0 \right\}$, so that the pole condition
requires $U$ to be analytic at $-i \omega$, thus
removing this pole from $U$ and the corresponding mode from the
solution. Note that in this simple one-dimensional example, the pole
condition can also be enforced by requiring the numerator of the
right term in~\eqref{eq:Udecomp} to vanish, leading to
\begin{equation}
  i \omega \hat{u}_{0} + \hat{u}_{0}' = 0 \ \leftrightarrow \ \partial_{t} u_{0} - u_{0}' = 0,
\end{equation}
which is the well known transparent boundary condition for the
one-dimensional wave equation, see for example~\cite{EngquistMajda}.

For more complex problems, an explicit decomposition of $U$
like~\eqref{eq:Udecomp} is usually not available.  However, the
Laplace transform can often still be decomposed into incoming and
outgoing parts in terms of path integrals in the complex plane,
see~\cite{RuprechtEtAl2008}, hence still allowing to define a region
$\mathbb{C}_{\rm in}$ and a pole condition based transparent boundary
condition. In particular, this is possible for the different types of
equations listed in Table~\ref{tab:Cin}.

\section{Discretization}\label{sec:Disc}
This section presents the employed
discretizations. Subsection~\ref{subsec:spaceDisc} describes the
discretization of the exterior domain with special semi-infinite
elements and how the pole condition is incorporated. The mapping
between the exterior elements and the corresponding reference element
introduces a generalized distance coordinate, along which the pole
condition is enforced. As the exterior elements are semi-infinite,
integrals arise that have one limit infinite. A mapping is
introduced, converting these integrals into proper integrals in the
Hardy space on the complex unit-disc. The discretization of the
interior uses standard finite elements and is not elaborated
further. Subsection~\ref{subsec:timeDisc} describes the employed
integration schemes, including the choices of the parameter of the
mappings.

\subsection{Space discretization}\label{subsec:spaceDisc}
The discretization of the exterior domain uses semi-infinite
trapezoids as proposed in~\cite{ZschiedrichPML,Zschiedrich03}.  The
construction of these exterior meshes is discussed
in~\cite{Kettner11}.  Possible other choices would be semi-infinite
triangles and rectangles as in~\cite{NannenSchaedle:11}. In any case
the exterior discretization has to be such that there is a uniform
distance variable. For the sake of simplicity, we assume that the
computational domain $\Omega$ is convex, although a generalization to
star-shaped domains should be possible.

\subsubsection{Transformation of exterior elements}
\begin{figure}[t]
  \centering
 \includegraphics[width=0.35\textwidth]{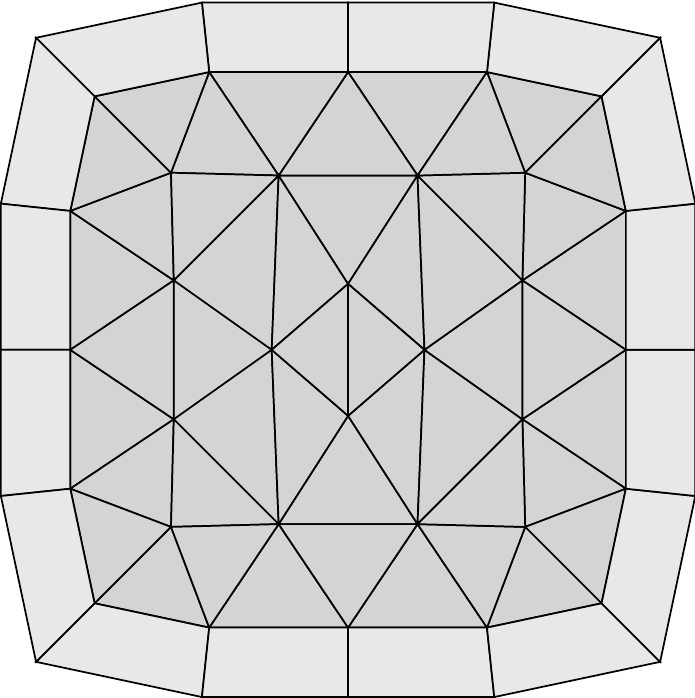}
  \caption{Employed basic mesh. Elements triangulating
    the interior domain are marked in dark grey while the trapezoids
    decomposing the exterior are marked in light grey. Finer meshes are
    generated by uniform refinements of the
    triangles and adding additional rays and
    exterior elements when new nodes on the boundary
    emerge.\label{fig:mesh} }
\end{figure}

\begin{figure}
  \centering
	\includegraphics[width=0.6\textwidth]{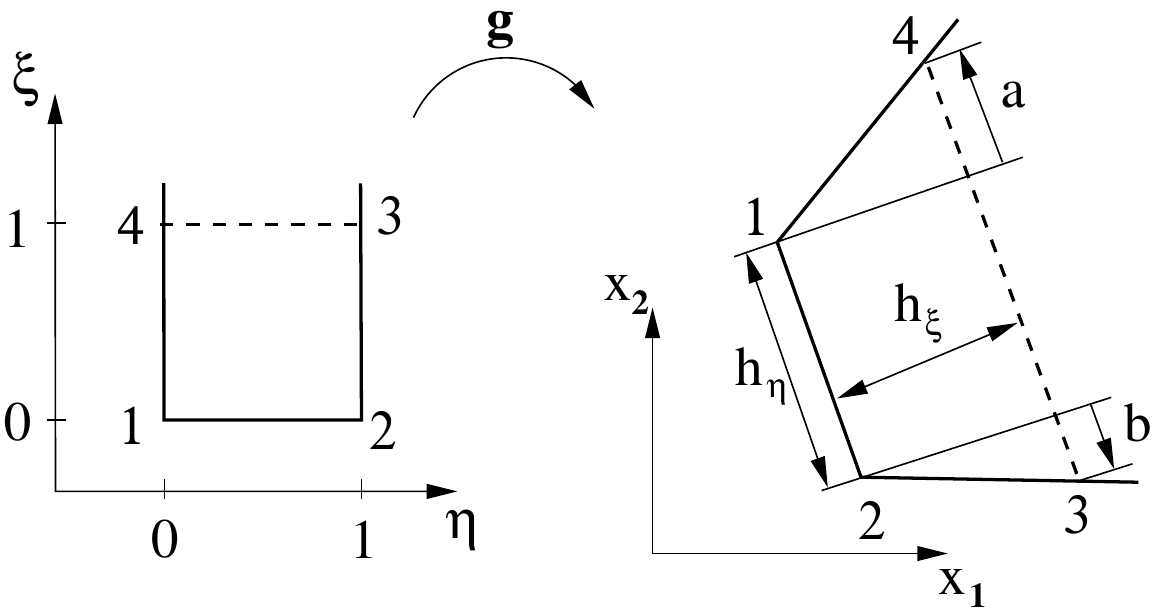}
	\caption{Transformation mapping the reference semi-infinite rectangle to an exterior trapezoidal element.\label{fig:trafo}}
\end{figure}

We mainly use the notation of~\cite{NannenSchaedle:11}. Figure
\ref{fig:mesh} sketches the used mesh including the exterior elements
and the mapping between the semi-infinite trapezoids in the exterior
and the corresponding reference element. Each of the semi-infinite
elements $T$ is the image of the semi-infinite reference rectangle
$[0,1]\times[0,\infty)$ under the bilinear mapping $\Tt{g}$ for a set
of parameters $(h_{\eta}, h_{\xi}, a,b)$. Denoting the Jacobian matrix
of $\Tt{g}$ by $J$ and its determinant by $|J|$, the mass and
stiffness integrals of the variational formulation, which may by found
in~\cite{NannenSchaedle:11}, are repeated below for the convenience of
the reader. Additionally the drift term is given.
\begin{equation}
  \label{eq:lokalInt}
  \begin{aligned}
  \int_{T} \nabla_{\Tt{x}} u \cdot \nabla_{\Tt{x}} v \,d\Tt{x}  &= 
  \int_{[0,1]\times[0,\infty]}
  J^{-T} \nabla_{\eta\xi} \tilde{u} \cdot J^{-T} \nabla_{\eta\xi}
  \tilde{v} |J| \, d (\eta,\xi),
  \\
  \int_{T} \Tt{d} \cdot \nabla_{\Tt{x}} u ~ v \, d\Tt{x} &= 
  \int_{[0,1]\times[0,\infty]} \Tt{d} \cdot J^{-T} \nabla_{\eta\xi}
  \tilde{u} ~ \tilde{v} |J| \, d(\eta,\xi).
  \\
  \int_{T} u ~ v \, d\Tt{x}  &=   \int_{[0,1]\times[0,\infty]}
  \tilde{u} ~ \tilde{v} |J| \, d (\eta,\xi),
\end{aligned}
\end{equation}
The Jacobian of the bilinear mapping $\Tt{g}$ and its determinant are
\begin{equation}
  J= R\left(
    \begin{array}{cc}
      h_\eta+(a+b)\xi& -b+(a+b)\eta 
      \\
       0& h_\xi
    \end{array}\right), \qquad
  |J| = h_{\xi}(h_{\eta}+(a+b)\xi), 
\label{eq:Jacobian}
\end{equation}
where $R$ is a rotation, $h_{\eta}$ is the
width of the trapezoid, $h_\xi$ is a scaling factor measuring the
distance to the boundary and $a$ and $b$ are signed distance
variables, compare for Figure~\ref{fig:trafo}.
\par
As in~\cite{NannenSchaedle:11} we use a
product ansatz $\tilde u(\eta,\xi) := u(\Tt{g}(\eta,\xi)) = \tilde u_\xi(\xi) \tilde u_\eta(\eta)$ 
on the reference strip, where $\tilde u_\xi$ and $\tilde u_\eta$ are functions
in $\xi$ and $\eta$, respectively. The stiffness term, for example, is given by
\begin{equation}
  \label{eq:stiffness}
  \begin{aligned}
  &  \int_{[0,1]\times[0,\infty]}
  J^{-T} \nabla_{\eta\xi} \tilde{u} \cdot J^{-T} \nabla_{\eta\xi}
  \tilde{v} |J| \, d (\eta,\xi) = 
  \\
  & \quad \int_0^1 \int_0^\infty 
  \left(
    \begin{array}{c}
      \partial_\eta\tilde{u}_\eta \tilde{u}_\xi\\
      \tilde{u}_\eta \partial_\xi\tilde{u}_\xi
    \end{array}
  \right)^T
  \left(
    \begin{array}{c c}
    \frac{ h_\xi^2+(b-(a+b)\eta)^2}{h_\eta + \xi(a+b)} & \frac{b-(a+b)\eta}{h_\xi}  \\
    \frac{b-(a+b)\eta}{h_\xi}  & \frac{h_\eta+(a+b)\xi}{h_\xi} 
  \end{array}
 \right)
  \left(
    \begin{array}{c}
      \partial_\eta\tilde{v}_\eta \tilde{v}_\xi\\
      \tilde{v}_\eta \partial_\xi\tilde{v}_\xi
    \end{array}
  \right) d\xi d\eta.
\end{aligned}
\end{equation}
Functions $\tilde{u}_\eta$ will be approximated using standard finite element basis functions $\{\phi_j(\eta)\}$, 
$j=1,\dots,N_\eta$ and therefor integrals over the radial coordinate $\eta$ can be evaluated by quadrature formulas. 

\subsubsection{Hardy Space}
\begin{figure}[t]
  \centering
  \includegraphics[width=0.4\textwidth]{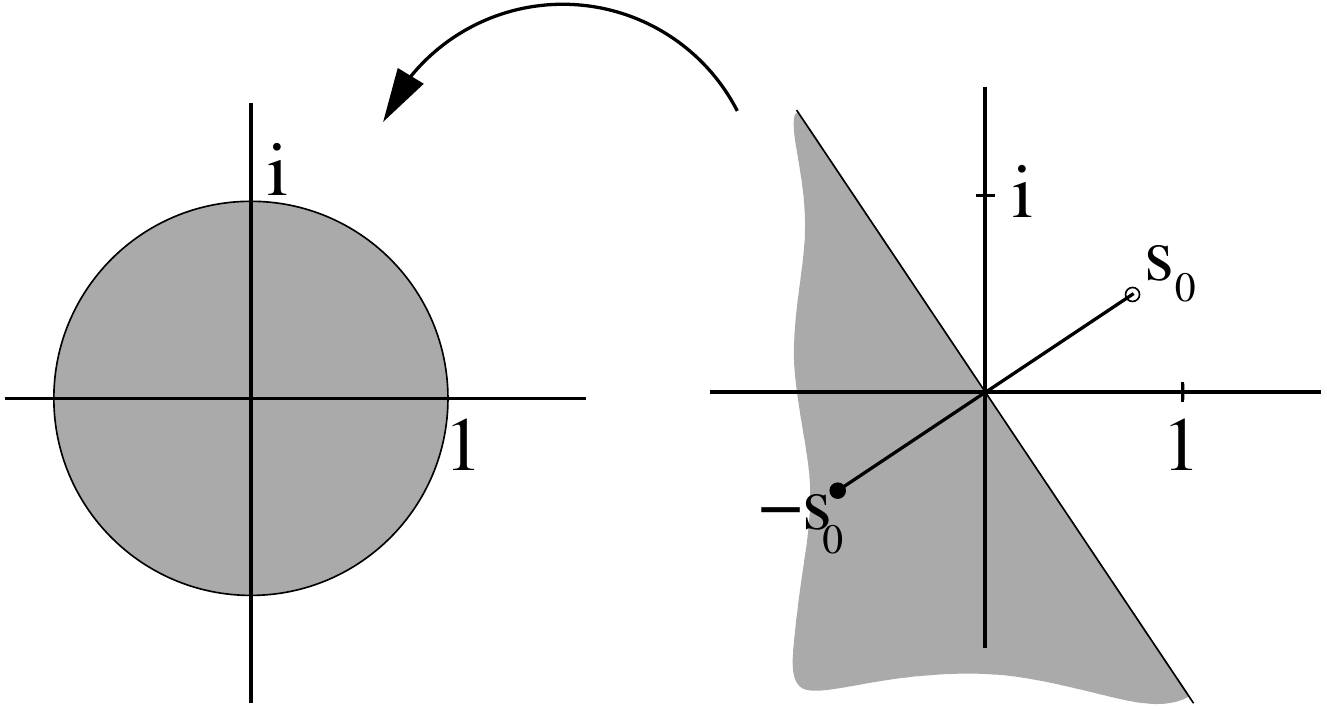}
  \caption{Sketch of the M\"obius transform $\mathcal{M} : \ H^-(P_{s_0}) \to H^+(D_{0})$, mapping functions defined on the complex half-plane $P_{s_{0}}$
  to functions defined on the complex unit-disc $D_{0}$.    
    The parameter $s_{0}$ has to be chosen
    such that $P_{s_{0}}$ corresponds to the $\mathbb{C}_{\rm in}$ suitable for the problem at hand.\label{fig:moebius}}
\end{figure}

Infinite integrals over the radial coordinate $\xi$ are transformed to finite integrals in the Hardy space
$H^+(D_{0})$ using the identity
\begin{equation}
  \label{eq:basisid} 
  \int_0^\infty \tilde{f}(\xi) \tilde{g}(\xi) d\xi = 
  -2 s_0 \frac{1}{2\pi} \int_{\partial D_{0}}
  (\mathcal{M}\mathcal{L}f)(\bar{z}) (\mathcal{M}\mathcal{L}g(z)) |dz|
\end{equation}
where $\mathcal{M}$ denotes the modified M\"obius transform
\begin{equation} 
  H^-(P_{s_0}) \rightarrow H^+(D_{0}) \, : \, F \mapsto \mathcal{M}F \mbox{ defined by } (\mathcal{M} F)(z):= F\left( s_0 \frac{z+1}{z-1}\right) \frac{1}{z-1}.
\end{equation} 
Here, $P_{s_{0}}$ denotes a half-plane in the complex plane, depending
on the parameter $s_{0}$, and $D_{0}$ denotes the complex unit-disc,
see Figure \ref{fig:moebius}. Further, $\mathcal{M}$ is an isomorphism
between the Hardy spaces $H^{-}(P_{s_{0}})$ and
$H^{+}(D_{0})$. Details can be found
in~\cite{Nannen2008,NannenSchaedle:11}.

The parameter $s_{0}$ has to be chosen such that the half-plane
$P_{s_{0}}$ coincides with the half-plane $\mathbb{C}_{\rm in}$ of
non-physical poles for the considered problem. As functions in the
space $H^{-}(P_{s_{0}})$ are analytic on the half-plane $P_{s_{0}}$
and $\mathcal{M}$ is an isomorphism, a function $\mathcal{M}F$ is
analytic on $D_{0}$ if and only if $F$ is analytic on $P_{s_{0}}$.

Hence the pole condition, stating that the Laplace transform
$\mathcal{L}(f)$ of some function $f$ has to be analytic on
$\mathbb{C}_{\rm in}$, is equivalent to the condition that
$\mathcal{M}\mathcal{L}(f)$ is analytic on $D_{0}$ for the correct
choice of the parameter $s_{0}$. In short, we established the
following sequence of reformulations of the pole condition
\begin{equation}
  \begin{aligned}
    f \ \textrm{satisfies pole condition} & :\Leftrightarrow F \ \textrm{has analytic extension to } \mathbb{C}_{\rm in} \\
    & \Leftrightarrow F \in H^{-}(P_{s_0}) \ \textrm{for correct choice of } s_{0} \\
    & \Leftrightarrow \mathcal{M} F \in H^{+}(D_{0}),
  \end{aligned}
\end{equation}
see~\cite{Nannen2008} for details.

In order to derive a formulation, which is easy to implement, some more transformations are required: Given a function $f$,  its image under $\mathcal{L}$ and $\mathcal{M}$
is decomposed into
\begin{equation} 
  \label{eq:decomUV}
  \mathcal{M} \mathcal{L} (f)(z) =
  \frac{1}{2s_0}( f_0 + (z-1)F(z)) =:
  \frac{1}{s_0}\mathcal{T}_- \left(
    \begin{array}{c} f_0 \\ F \end{array}\right)(z) 
\end{equation}
in order to get a local coupling with the boundary data $f_0$. As the Laplace transform maps differentiation to multiplication by $s_0 (z+1)/(z-1)$, straightforward calculation yields
\begin{equation}
  \label{eq:Tp}
  \mathcal{M} \mathcal{L} (f')(z) = \frac{1}{2}( f_0 + (z+1)F(z)) =: \mathcal{T}_+ \left(
    \begin{array}{c} f_0 \\ F \end{array}\right)(z).
\end{equation}
For motivating this coupling, note that it follows from general theory
on Laplace transforms that if $F$ is the Laplace transform of $f$ one
has $\lim_{s\to \infty} sF(s) = f(0)$, whenever $f(0)$ exists. For
details on the decomposition we again refer
to~\cite{Nannen2008,NannenSchaedle:11}.

It remains to take care of terms which contain multiplications by
$\xi$ and $((a+b)\xi+h_\eta)^{-1}$ in~\eqref{eq:stiffness}. To this
end, an additional operator $\mathcal{P} : H^+(D_{0}) \to H^+(D_{0})$
(multiplication by $\xi$) \footnote{This operator is denoted by $D$
in~\cite{NannenSchaedle:11}.}  is implicitly defined by
\begin{equation*} 
  \mathcal{M} \mathcal{L} \{ (\cdot ) f(\cdot) \}=\mathcal{M} \left\{
    -\left( \mathcal{L} f(\cdot)\right)' \right\}
  = s_0^{-1}\mathcal{P} \left(\mathcal{M} \mathcal{L} f\right).
\end{equation*}
Direct calculations yield
\begin{equation}
  \label{eq:defDop} \left(\mathcal{P} F \right)(z) =
  \frac{(z-1)^2}{2} F'(z)+\frac{z-1}{2} F(z),\qquad F \in
  H^+(D_{0}).
\end{equation}
Assembling the discrete system involves basically the assembly of discrete counterparts of the operators
$\mathcal{T}_{+}$, $\mathcal{T}_{-}$ and $\mathcal{P}$.

\subsubsection{Choosing Basis Functions}
Up to now, all transformations were on a continuous level and no
approximations were made. Because functions that are analytic on the complex
unit-disc can be expanded in power series, the set of monomials
$\left\{ z^{j} \right\}_{j=0}^{\infty}$ constitutes a basis of
$H^{+}(D_{0})$. By using the space spanned by a finite number of
monomials $\{z^j\}_{j=0}^{\Nxi}$ as test and ansatz space, one obtains
finite dimensional approximations of $\mathcal{T}_{\pm}$
\begin{align}
  \label{eq:Mat}
  T_{\Nxi,+} =
  \left(\begin{array}{cccc}
    1 & 1 & &   \\
      & 1 & 1 &  \\
      & & \ddots & \ddots  \\
      & & & 1
  \end{array}\right), \ 
  T_{\Nxi,-} =
  \left(\begin{array}{cccc}
    1 & -1 & &   \\
      & 1 & -1 &  \\
      & & \ddots& \ddots \\
      & & & 1
  \end{array}\right),
\end{align}
and $\mathcal{P}$
\begin{align}
  P_{\Nxi} =
  \left(\begin{array}{cccc}
    -1 & 1 & &   \\
    1  & -3 & 2 &  \\
      & \ddots& \ddots& \ddots  \\
      & & (N_{\xi}-1) & -(2N_{\xi}+1)
  \end{array}\right).
\end{align}
The Hardy space monomials can also be transformed back to give a representation of the corresponding ansatz and test functions in physical space,
see~\cite{Nannen2008}. To define the local stiffness matrix, set for the $\xi$-integrals
\begin{equation}
  \begin{aligned}
    L_{\xi,11}^{(-1)}&:=-2 ~ T_{\Nxi,-}^\top \left(h_\eta s_0~I+(a+b) P_{\Nxi}\right)^{-1} T_{\Nxi,-} ,\\
    L_{\xi,12}^{(0)} &:=-2 ~ T_{\Nxi,-}^\top  T_{\Nxi,+} , \quad
    L_{\xi,21}^{(0)} :=-2 ~ T_{\Nxi,+}^\top  T_{\Nxi,-} ,\\
    L_{\xi,22}^{(1)} &:= -2h_\eta T_{\Nxi,+}^\top T_{\Nxi,+} , \quad
    L_{\xi,22}^{(0)}  :=-2 (a+b) T_{\Nxi,+}^\top P_{\Nxi} T_{\Nxi,+}. 
  \end{aligned}
  \label{eq:Ldefs}
\end{equation}
Here the superscript counts the leading order in $s_0$ and the subscripts correspond to
the position in the matrix in~\eqref{eq:stiffness}. For the $\eta$-integrals set
\begin{equation}
  \begin{aligned}
    L_{\eta,11} &:= \Big(\int_0^1 \phi'_i(\eta)\left(h_\xi + \frac{ ((a+b)\eta-b)^2}{h_\xi}\right) \phi'_j(\eta)\Big)_{i,j=1}^{N_\eta}, \\
    L_{\eta,12} &:= \Big(\int_0^1 \phi'_i(\eta) \frac{b-(a+b)\eta}{h_\xi} \phi_j(\eta)\Big)_{i,j=1}^{N_\eta}, \\
    L_{\eta,21} &:= \Big(\int_0^1 \phi_i(\eta) \frac{b-(a+b)\eta}{h_\xi}\phi'_j(\eta)\Big)_{i,j=1}^{N_\eta}, \\
    L_{\eta,22} &:= \Big(\int_0^1 \phi_i(\eta)  \frac{1}{h_\xi} \phi_j(\eta)\Big)_{i,j=1}^{N_\eta}.
  \end{aligned}
\end{equation}
For equations with second order temporal derivatives, the parameter $s_{0}$ is chosen to be 
frequency dependent in Fourier space, to be precise $s_{0} = i\omega$,
translating back to $\partial_{t}$ in physical space. To avoid the inversion in $L_{\xi,11}^{(-1)}$ in~\eqref{eq:Ldefs}, additional unknowns are introduced such that the local stiffness matrices are given by
\begin{align}
  \label{eq:LocLaplaceInfQuad}
  \nonumber
  L^{(0)}_{loc} &= \left[
  \begin{array}{cc}
    L_{\eta,22} \otimes L_{\xi,22}^{(0)} + L_{\eta,12} \otimes L_{\xi,12}^{(0)} +  L_{\eta,21}\otimes L_{\xi,21}^{(0)} &
    -2 L_{\eta,11}  \otimes T_{\Nxi,-}^\top \\
    2 I \otimes T_{\Nxi,-} & -2 I \otimes (a+b) M_{\Nxi}
  \end{array}
  \right] \\
  L^{(1)}_{loc} &= \left[
  \begin{array}{cc}
     L_{\eta,22} \otimes L_{\xi,22}^{(1)} & 0 \\
     0 & -2 I \otimes h_\eta I
  \end{array}
  \right].
\end{align}
Similarly, local mass matrices corresponding to the mass integral
in~\eqref{eq:lokalInt} are given by
\begin{equation}
  M^{(-1)}_{loc} :=\left[\begin{array}{cc} 
      M_{\xi}^{(-1)} \otimes M_\eta & 0 \\
      0 & 0 \end{array}
  \right], \quad
  M^{(-2)}_{loc} :=\left[\begin{array}{cc} 
      M_{\xi}^{(-2)} \otimes M_\eta  & 0 \\
      0 & 0  \end{array}
  \right]\label{eq:Mloc}
\end{equation}
where
\begin{align}
  \begin{aligned}
    M_{\xi}^{(-1)}&:=-2 h_\xi h_\eta T_{\Nxi,-}^\top T_{\Nxi,-}
    \\
    M_{\xi}^{(-2)}&:=-2 h_\xi (a+b) T_{\Nxi,-}^\top P_{\Nxi}
    T_{\Nxi,-}
  \end{aligned}
\mbox{ and }
  M_\eta := \Big(\int_0^1 \phi_i(\eta) \phi_j(\eta)\Big)_{i,j=1}^{N_\eta}.
\end{align}
For the drift term set
\begin{equation}
  \begin{aligned}
    D_{\xi,1}^{(0)}  &:=  -2 T_{\Nxi,-}^\top T_{\Nxi,+}\\
    D_{\xi,2}^{(-1)} &:=  -2(a+b) T_{\Nxi,-}^\top P_{\Nxi}T_{\Nxi,+}\\
    D_{\xi,3}^{(-1)} &:=  -2 T_{\Nxi,-}^\top T_{\Nxi,-}\\
  \end{aligned}
\end{equation}
and
\begin{equation}
  \begin{aligned}
    D_{\eta,1} &:= h_\eta \tilde{d}_2\Big(\int_0^1 \phi_i(\eta)  \phi_j(\eta)\Big)_{i,j=1}^{N_\eta} \\
    D_{\eta,2} &:=  \tilde{d}_2\Big(\int_0^1 \phi_i(\eta)  \phi_j(\eta)\Big)_{i,j=1}^{N_\eta} \\
    D_{\eta,3} &:=  \Big(\int_0^1 \phi_i(\eta) 
    (\tilde{d}_2(b-(a+b)\eta)+\tilde{d}_1h_\xi) \phi'_j(\eta)\Big)_{i,j=1}^{N_\eta} \\
  \end{aligned}
\end{equation}
where $(\tilde d_1, \tilde d_2)^T = \Tt{\tilde{d}} = R\Tt{d}$ is the rotated $\Tt{d}$ vector. The local
drift matrices are then given by
\begin{equation}
D_{loc}^{(0)} :=\left[\begin{array}{cc} D_{\xi,1}^{(0)} \otimes
    D_{\eta,1} & 0 \\ 0 & 0 \end{array} \right],
D_{loc}^{(-1)}:=\left[\begin{array}{cc} D_{\xi,2}^{(-1)} \otimes
    D_{\eta,2} + D_{\xi,3}^{(-1)} \otimes D_{\eta,3} & 0 \\ 0 &
    0 \end{array} \right]
\end{equation}
In the computational domain $\Omega$ standard local finite element matrices
$M_{loc}^{(0)}$, $D_{loc}^{(0)}$ $L_{loc}^{(0)}$ without the $s_0$-parameter are obtained. By assembling the local matrices to global matrices, a spatial semi-discretization of~\eqref{eq:u} is obtained
\begin{equation}
  \label{eq:semidiscEq}
  \begin{aligned}
    p(\partial_t) &
    \left(M^{(0)} +\frac{1}{s_0} M^{(-1)} +\frac{1}{s^2_0} M^{(-2)}\right) u(t) 
    = \left(L^{(0)} + s_0 L^{(1)} \right) u(t) + \\&
    \left(D^{(0)} + \frac{1}{s_0} D^{(-1)}\right) u(t) 
    - k^2 \left(M^{(0)} + \frac{1}{s_0}  M^{(-1)} +\frac{1}{s^2_0} M^{(-2)}\right) u(t)
  \end{aligned}
\end{equation}
where $u(t)$ is the time-dependent vector of degrees-of-freedom, including the coefficients of the monomial basis functions of the subset of $H^{+}(D_{0})$
providing the boundary condition.

\subsection{Time discretization}\label{subsec:timeDisc}
All conducted simulations rely on the method-of-lines approach: The
PDE at hand is first discretized in space, as described in
subsection~\ref{subsec:spaceDisc}, leading to the
ODE~\eqref{eq:semidiscEq} for the coefficients. This equation is then
integrated in time using different time-stepping schemes indicated below.

\subsubsection{Schr\"odinger's equation}
For Schr\"odinger's equation,~\eqref{eq:semidiscEq} is solved by the
second order accurate, A-stable trapezoidal/mid-point rule. Denoting
by $u^n$ the approximation to $u(n\dt)$ at $t=n\dt$ for some time-step
size $\dt$, the discretization reads
\begin{equation}
\partial_t  M u(t) = -i c^2 L u(t) \leftrightarrow 
M\frac{u^{n+1}-u^{n}}{\dt} = -i c^2 L \frac{u^{n+1} + u^{n}}{2}
\end{equation}
where the mass and stiffness matrix are given by $M=M^{(0)} +s_0^{-1}
M^{(-1)} +s^{-2}_0 M^{(-2)}$ and $L=L^{(0)} + s_0 L^{(1)}$. As in the
one-dimensional case, non-physical solutions correspond to poles in
the first quadrant, hence $s_0$ is chosen in the third quadrant and
set to
\begin{equation}
  s_{0} = -1 - i,
\end{equation}
in order to exclude poles in the region $\mathbb{C}_{\rm in}$ indicated in Table \ref{tab:Cin}.

\subsubsection{Drift-Diffusion equation}
In the examples for the drift-diffusion equation,~\eqref{eq:semidiscEq} is integrated with the A-stable Radau IIA method with three stages of order five, see~\cite[Sec. IV.5]{HairerWannerII}. The parameter $s_0$ is chosen to be real and negative, such that the poles in the positive half-plane are excluded, corresponding to non-physical exponentially increasing solutions. We set
\begin{equation}
	s_{0} = -5,
\end{equation}
hence poles with positive real part, see Table~\ref{tab:Cin}, are excluded. The chosen value of $s_{0}$ produces good results, but some optimization is probably still possible. However,
the sensitivity of the results to the specific value is rather low, as long as the correct half-plane is excluded.

\subsubsection{Wave equation}
For the wave equation, poles in the lower complex half-plane have to be excluded, see~\ref{tab:Cin}.
For the wave equation it is $p(\partial t) = \partial_{tt}$, corresponding to $p(\omega) = -\omega^{2}$ 
in frequency space. As in the one-dimensional case, we choose the parameter $s_{0}$ to be frequency dependent, 
setting
\begin{equation}
  s_{0} = i \omega.
\end{equation}
Transforming back to physical space yields
\begin{equation}
  \partial_{tt} M^{(0)} u  - \partial_{t} M^{(-1)}u + M^{(-2)} u  = 
    L^{(0)} u  - \partial_t L^{(1)} u.
\end{equation}
Discretization is done again with the implicit trapezoidal rule, resulting in
\begin{equation}
  \begin{aligned}
    M^{(0)} \frac{u^{n+1}-2u^n+u^{n-1}}{\dt^2} - &M^{(-1)}\frac{u^{n+1}
      - u^{n-1}}{2\dt} + M^{(-2)} \frac{u^{n+1}+2u^n+u^{n-1}}{4} = 
    \\
   & L^{(0)}\frac{u^{n+1}+2u^n+u^{n-1}}{4} - L^{(1)} \frac{u^{n+1} - u^{n-1}}{2\dt}. 
  \end{aligned}
  \label{eq:trapRule}
\end{equation}

\section{Numerical results}\label{sec:numExam}
The computational domain in all simulations is a square $[-4, 4]
\times [-4, 4]$ in the two-dimensional plane with slightly smoothed
corners, see Figure~\ref{fig:mesh}. Sketched in light gray are the
trapezoidal elements spanned by the rays in the exterior domain while
the darker triangles correspond to the triangulation of the interior
domain. In order to obtain higher resolutions, the shown mesh is
refined using up to five uniform refinement steps. As the original
mesh is very coarse, no errors are reported for simulations on the
unrefined grid, because at least lower order finite elements do not
produce reasonable solutions there.

\subsection{Schr\"odinger's equation}
For $p(\partial_{t}) = i \partial_{t}$, $c = k = 0$ and $\Tt{d} = (0,0)^T$, 
equation~\eqref{eq:u} yields Schr\"odinger's equation. For this case,
exact solutions of the form
\begin{equation}
	u_{\alpha}(x, y, t; \alpha) = \frac{i}{4t + i}  \exp\left( \frac{-i \left( x^{2} + y^{2}  \right) 
          - \alpha \left( x + y \right) - 2 \alpha^{2} t  }{4t + i} \right)
\end{equation}
with a parameter $\alpha$ are available. We use a superposition of two such solutions, that is
\begin{equation}
	u(x,y,t) = u_{\alpha}(x,y,t; \alpha=1.4) + u_{\alpha}(x,y,t; \alpha=-2),
\end{equation}
and employ $u(x,y,0)$ as initial value. The discretization in space employs finite elements of orders one to four in the interior.
Simulations are run until $T = 10$ with time-steps $\Delta t = 1/800, 1/1600, 1/3200, 1/6000$ on  meshes refined up to five times and
for values of $\Nxi$ (coefficients per ray) between $\Nxi = 1$ and $\Nxi = 20$. Output is generated at two hundred points in time, distributed equally over the time interval
$[0, 10]$.
\begin{figure}[t]
  \centering
  \includegraphics[scale=1]{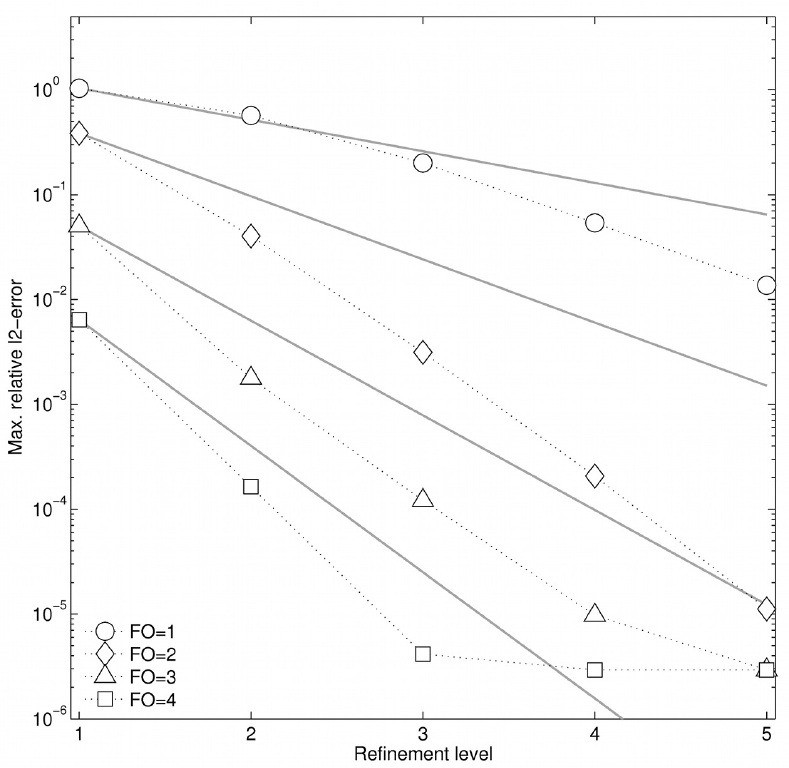}
  \caption{Verification of the spatial order of convergence for the Schr\"odinger equation. Shown is the maximum of the relative $l_2$-errors over all outputs versus the refinement level of the mesh for finite elements of order one to four. The employed time-step is $\Delta t = 1/6000$ and $\Nxi=20$ coefficients along each ray are used. As a guide to the eye, lines with slopes one to four have been added. \label{fig:ScEqSpatConv}}
\end{figure}
\par
Figure~\ref{fig:ScEqSpatConv} shows the maximum relative $l_2$-error over all generated outputs versus the refinement levels of 
the mesh for finite elements of order one to four. All elements converge with the expected rate or 
better until the error saturates at about $3 \times 10^{-6}$ in the case of the higher order finite elements. At this point, the error from the
temporal discretization starts dominating, compare for Table~\ref{tab:ScEq-TempConv}, and increasing the accuracy of the spatial discretization
yields no more improvement unless the accuracy of the time-discretization is also increased.
\begin{table}[t]
	\centering
	\begin{tabular}{|c|c|c|} \hline
		time-step & log10(error) &  Conv. Rate \\ \hline
		1/800 & -3.8 & -- \\
		1/1600 & -4.4 & 2.0 \\
		1/3200 & -5.0 & 2.0 \\
		1/6000 & -5.5 & 2.0 \\ \hline
	\end{tabular}
	\caption{Maximum relative $l_2$-error over all generated outputs for the Schr\"odinger equation, depending on the time-step size. The simulation used finite elements of order four, a five times refined mesh and $\Nxi = 20$ coefficients along each exterior ray.\label{tab:ScEq-TempConv}}
\end{table}
\par 
Table \ref{tab:ScEq-TempConv} shows the maximum relative $l_2$-error versus the length of the time-step, confirming the convergence rate of two
expected from the employed second order accurate trapezoidal rule. 
\begin{figure}
  \centering
  \includegraphics[scale=1]{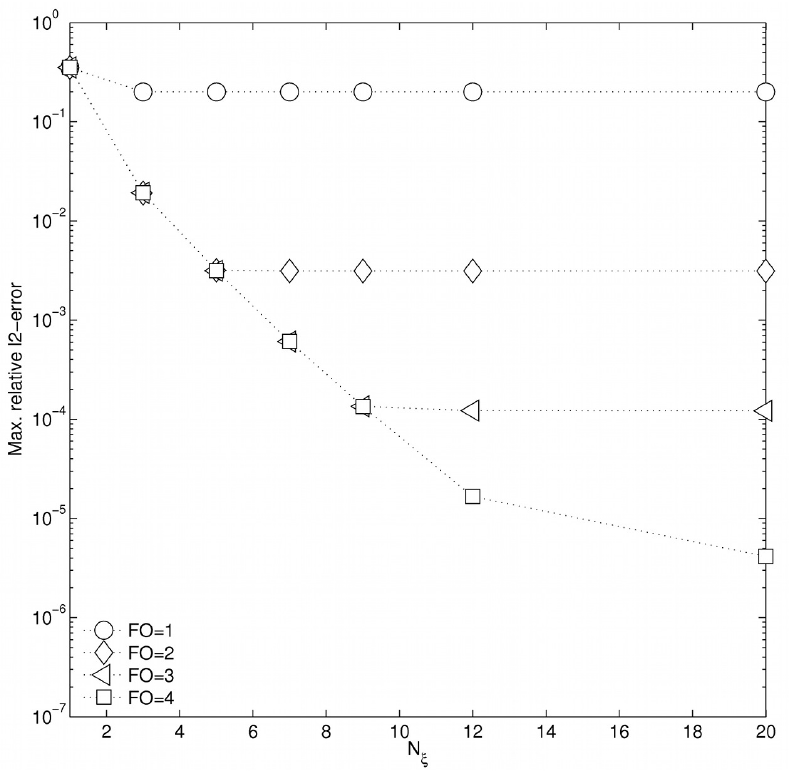}
  \caption{Relative $l_2$-error depending on the number of coefficients $\Nxi$ per ray in the exterior domain for Schr\"odinger's equation. Maximum error over all generated outputs for finite elements of order one to four on a three times refined mesh.\label{fig:ScEq-ErrorDec}}
\end{figure}
\begin{figure}
	\centering
	 \includegraphics[scale=1]{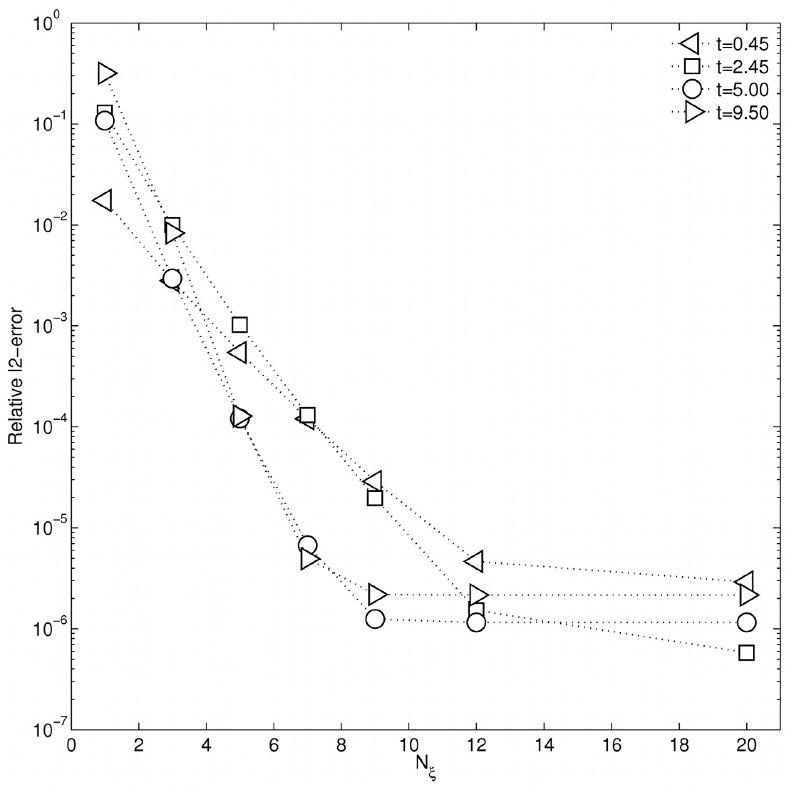}
	\caption{Relative $l_2$-error at four fixed points in time depending on the number of coefficients $\Nxi$ per ray for the simulation with finite elements of order four. In all cases, a time-step of $\Delta = 1/6000$ has been used.\label{fig:ScEq-ErrorDecDiffT}}
\end{figure}
\par
Figures~\ref{fig:ScEq-ErrorDec} and \ref{fig:ScEq-ErrorDecDiffT} show how the error decays with increasing $\Nxi$. The former shows the maximum relative $l_2$-error for simulations with finite elements of order one to four, a three times refined mesh and a time-step of $\Delta = 1/6000$. In all cases, the error decays super-algebraically
with $\Nxi$ until it saturates at the level of the respective spatial or temporal discretization error. Figure~\ref{fig:ScEq-ErrorDecDiffT} shows the error at four different points in time for the fourth order elements. The error decays super-algebraically at all four points in time, even at later times where most of the solution has left the domain.
\begin{figure}[t]
  \centering
  \includegraphics[scale=1]{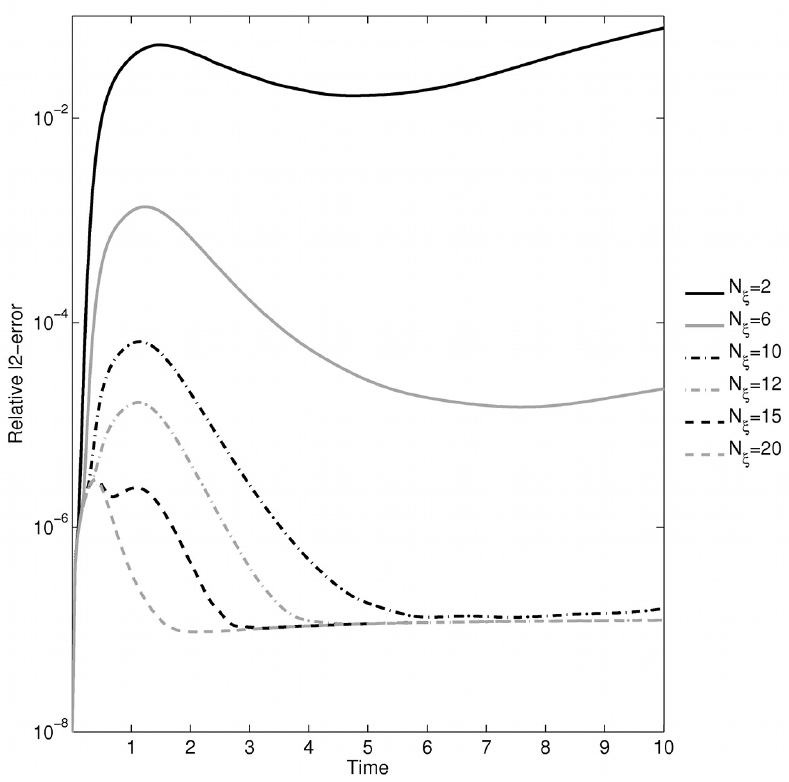}
  \caption{Relative $l_2$-error over time for the Schr\"odinger equation for six different values of $\Nxi$. The simulation used finite elements of order four, a five times refined mesh and a time-step $\Delta t = 1/6000$. \label{fig:ScEq-ErrorVsTime}}
\end{figure}
\par
Figure~\ref{fig:ScEq-ErrorVsTime} shows the relative $l_2$-error over time for different values of $\Nxi$. In all cases, the error 
increases as the wave packets hit the boundary of the computational domain and subsequently decays to a level determined by the number of 
coefficients per ray $\Nxi$. The error decays faster for larger values of $\Nxi$, but for values of $\Nxi = 10$ or larger, the 
levels at which the error saturates and in particular the error at the end of the simulation is identical.

%
% Drift-diffusion equation
%
\subsection{Drift-diffusion equation}
Setting $p(\partial_{t}) = \partial_{t}$ and $k=0$ in~\eqref{eq:u} yields the drift-diffusion equation. Note that the heat 
equation is included here as the special case $\Tt{d}=(0,0)^{\rm T}$. An analytic solution is given by
\begin{equation}
	u(x,y,t) = \frac{1}{t} \exp\left(  - \frac{1}{4 t c^{2}} \left[  \left(x - d_{1} t  \right) ^{2} 
          + \left(   y - d_{2} t  \right)^{2}\right] \right).
\end{equation}
Set
\begin{equation}
\Tt{d} = \left( d_{1} , d_{2} \right) = \left( 1.5 , 1.5 \right), \quad \text{and} \ c = 0.5
\end{equation}
and start the integration at $t_{0} = 0.2$ with initial value $u(x,y,t_{0})$. This yields a Gaussian function with a peak 
initially close to the origin which is subsequently advected to the upper right corner of the square while being spread 
out by diffusion. Integration in time is done by the fifth order Radau IIA(5) scheme. The simulations are run until $T=5$ with
finite elements of order one to six, on meshes refined up to four times, values of $\Nxi$ between $1$ and $51$
and time-steps ranging from $\Delta = 1/10$ to $\Delta t = 1/160$.
\begin{figure}
  \centering
  \includegraphics[scale=1]{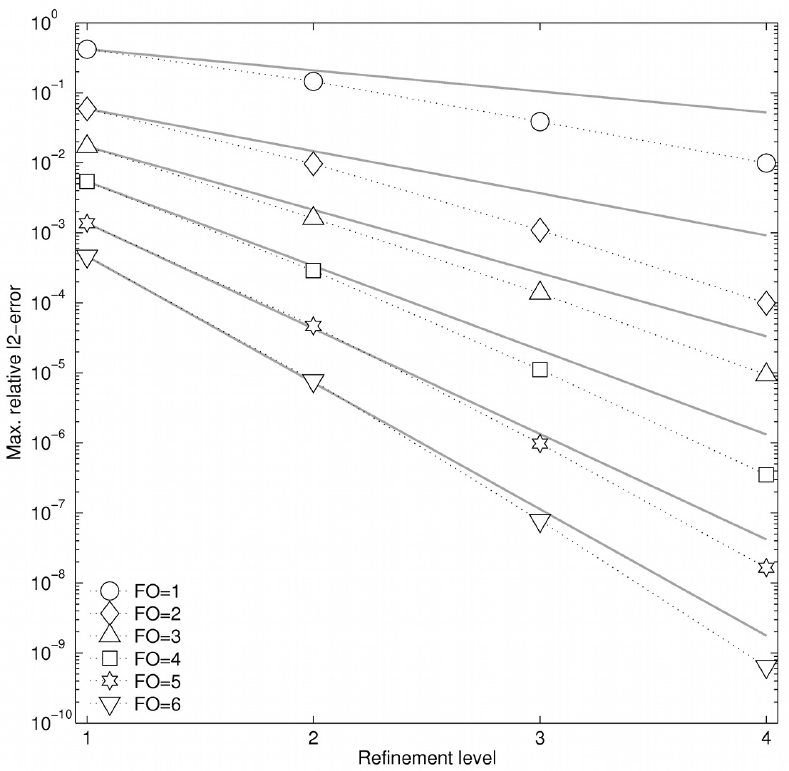}
    \caption{Verification of the spatial order of convergence for the drift-diffusion equation integrated until $T=5$. Shown is the maximum of the relative $l_2$-errors over all generated outputs. The number of coefficients per ray is $\Nxi = 51$ and the time-step is $\Delta t = 1 / 160$.\label{fig:DdEq-SpatConv}} 
 \end{figure}
\par Figure~\ref{fig:DdEq-SpatConv} shows the maximum relative $l_2$-error over all generated outputs versus the refinement level of the mesh. The
number of coefficients per ray is $\Nxi = 51$ and the time-step is $\Delta t = 1/160$. As a guide to the eye, lines with slopes from one 
to six are added. All elements converge with the expected rate or better, confirming again that the pole condition does not compromise the 
order of convergence of the spatial discretization.
\begin{table}[t]
	\centering
	\begin{tabular}{|c|c|c|} \hline
		time-step & log10(error) &  Conv. Rate \\ \hline
		1/10 & -3.6 & -- \\
		1/20 & -5.0 & 4.6 \\
		1/30 & -5.9 & 4.9 \\
		1/40 & -6.5 & 4.9 \\
		1/60 & -7.4 & 4.9 \\
		1/80 & -8.0 & 5.0 \\
		1/160& -9.2 & 4.0 \\ \hline
	\end{tabular}
	\caption{Maximum relative $l_2$-error over all generated outputs depending on the time-step size for finite elements of order six, a four times refined
	mesh and $\Nxi = 51$ coefficients per ray.}
	\label{tab:DdEq-TempConv}
\end{table}
\par 
Table~\ref{tab:DdEq-TempConv} shows the maximum relative $l_2$-error versus the length of the employed time-step for finite elements of order six,
a four times refined mesh and $\Nxi = 51$. From $\Delta t = 1/10$ to $\Delta t = 1/20$, a slightly reduced convergence rate is observed, probably because the time-step size is still in the pre-asymptotic regime. The reduced convergence rate in the last refinement is because the error approaches the spatial discretization error, see 
Figure~\ref{fig:DdEq-SpatConv}. Beside that, the expected fifth order convergence is observed, demonstrating that the pole condition can not only preserve the accuracy
of high order finite elements but also of high order integration schemes.
\begin{figure}
  \centering
  \includegraphics[scale=1]{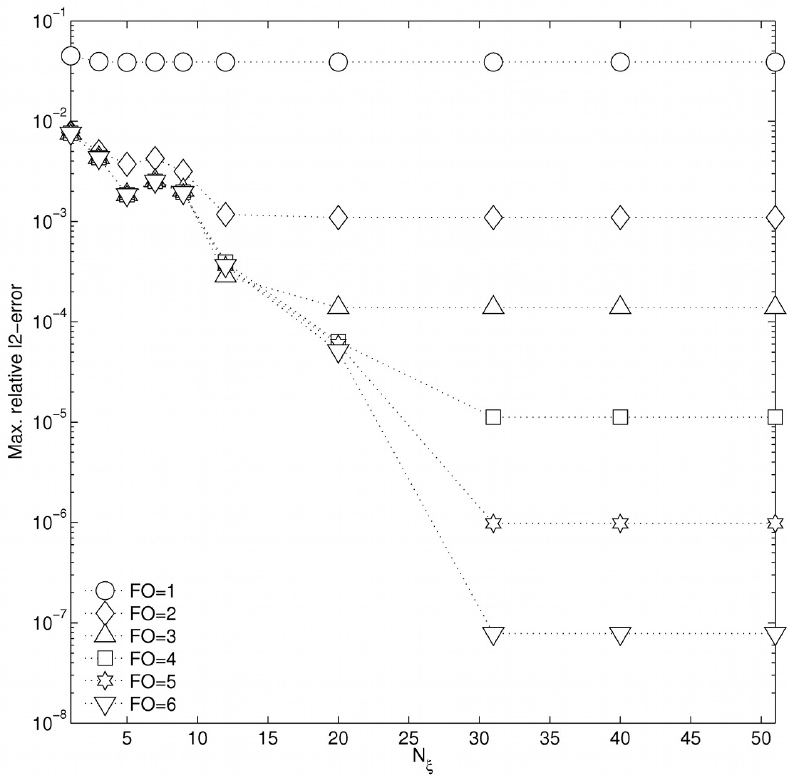}
    \caption{Relative $l_2$-error depending on the number of coefficients $\Nxi$ per ray in the exterior domain for the drift-diffusion equation. Maximum error over all generated outputs for finite elements of order one to six on a three times refined mesh.\label{fig:DdEq-ErrorDec}}
\end{figure}
\begin{figure}
	\centering
 	 \includegraphics[scale=1]{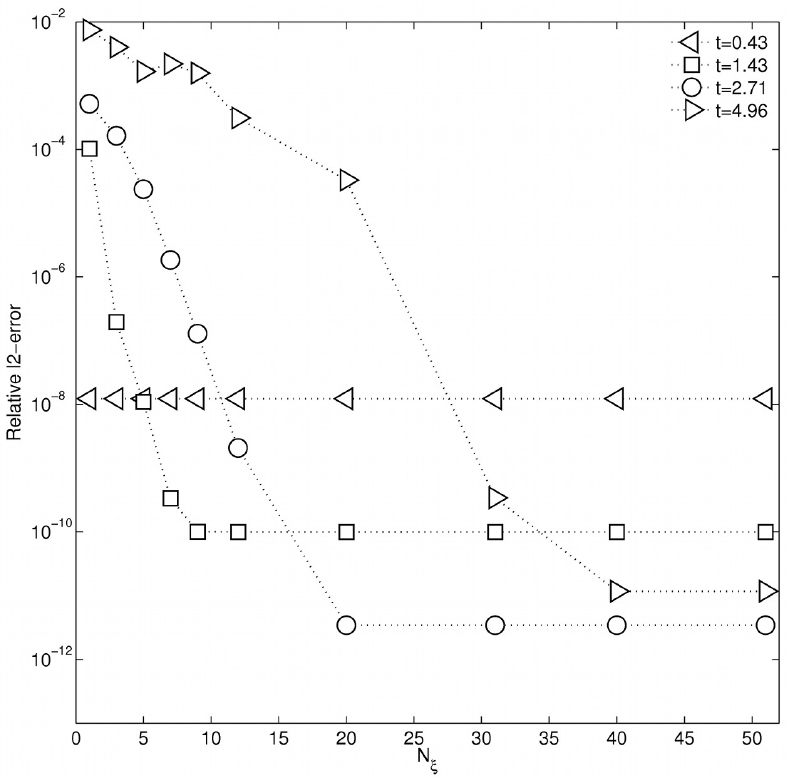}	
	\caption{Relative $l_2$-error depending on $\Nxi$ at four fixed points in time for the simulation with finite elements of order six. In all cases, a time-step $\Delta t = 1/160$ has been used.\label{fig:DdEq-ErrorDecDiffT}}
\end{figure}
\par 
Figures~\ref{fig:DdEq-ErrorDec} and~\ref{fig:DdEq-ErrorDecDiffT} show the relative $l_2$-error versus the number of coefficients $\Nxi$. The former shows the maximum error
over all outputs for finite elements of order one to six, a three times refined mesh and a time-step $\Delta t = 1/160$. In contrast to Schr\"odinger's equation,
for small values of $\Nxi$ there is only a minor decrease of the error. Also, below $\Nxi$, the error is not decreasing monotonically with the
number of coefficients. After $\Nxi= 10$, rapid super-algebraically decrease of the error is again observed until the error saturates at a level
determined by the accuracy of the spatial discretization. Note that $\Nxi = 30$ coefficients per ray are sufficient here to reduce the boundary condition error
to the level of the spatial discretization error in all cases. Figure~\ref{fig:DdEq-ErrorDecDiffT} shows the relative error at four different points in time. Again the error generally decays
super-algebraically with the number of coefficients, but now the decay rates are noticeably lower at later points in time. Note that at $t=0.43$, the Gauss peak has
not yet reached the boundary, so that the boundary condition has no visible effect on the error at this time.
\begin{figure}
  \centering
  \includegraphics[scale=1]{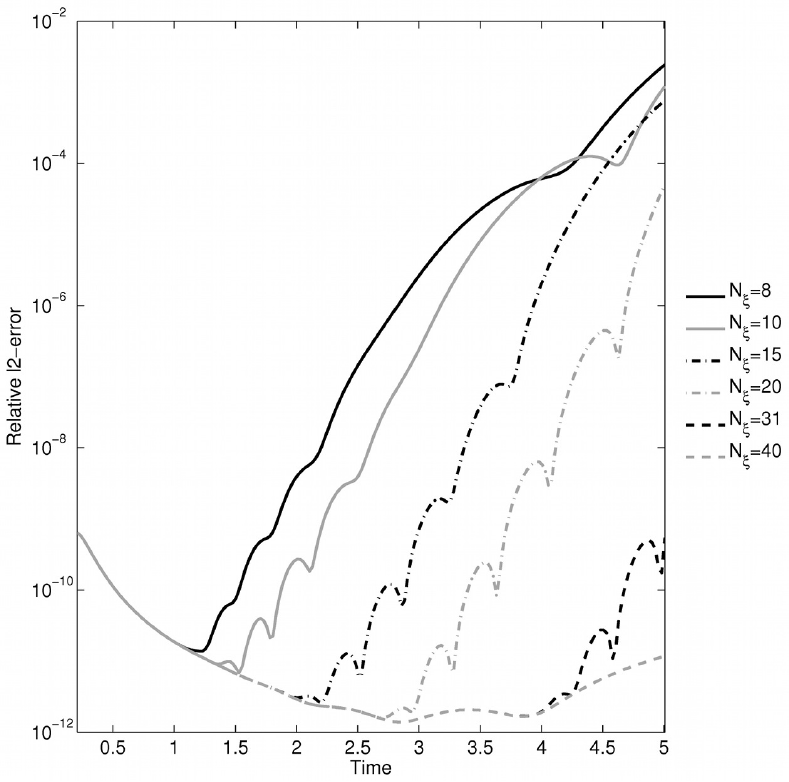}  
  \caption{Relative $l_2$-error over time for the drift-diffusion equation for different values of $\Nxi$ at a four times refined mesh, finite elements of order six and a time-step $\Delta t = 1/160$.\label{fig:DdEq-ErrorVsTime}}
\end{figure}
\par Figure \ref{fig:DdEq-ErrorVsTime}  shows the relative $l_2$-error over time for different values of $\Nxi$, finite elements of 
order six and a four times refined mesh. As for Schr\"odinger's equation, the error starts increasing at some point in time, but the increase starts later as the number of
coefficients increases. On the other hand, for the simulations with $\Nxi = 15$ or less coefficients, the errors at the end of the simulation are about the same and
only for larger values of $\Nxi$ a significantly reduced error is observed at the end of the simulation. Together with Figure~\ref{fig:DdEq-ErrorDec}, this illustrates that
for the drift-diffusion example, the error is not monotonically decreasing with $\Nxi$ for small values of $\Nxi$ and a certain minimum number of coefficients per ray is required before the onset of the super-algebraic decay.

\subsection{Wave and Klein-Gordon equation}
For $p(\partial_{t}) = \partial_{tt}$ and $\Tt{d}=(0,0)^{\rm T}$, \eqref{eq:u} becomes the Klein-Gordon equation, containing the wave 
equation as the special case $k=0$. While the pole condition could successfully provide TBCs for both equations in the one-dimensional case as well as in 
a two-dimensional wave-guide problem, stability problems arise in the fully two-dimensional case, rendering the pole condition in the here presented 
form inapplicable to both equations for finite elements of order two or higher. Resolving these issues is planned for future research.\par
Below, the instability is documented briefly. Use an initial distribution
\begin{equation}
	u(x, y) = \exp \left( -2x^{2} - 2 y^{2} \right),
\end{equation}
finite elements of order one to four and up to four refinement steps for the mesh. Integrate in time using implicit 
trapezoidal rule with a time-step $\Delta t = 1/1280$ until $T=100$.
\begin{figure}
\centering
\includegraphics[scale=1]{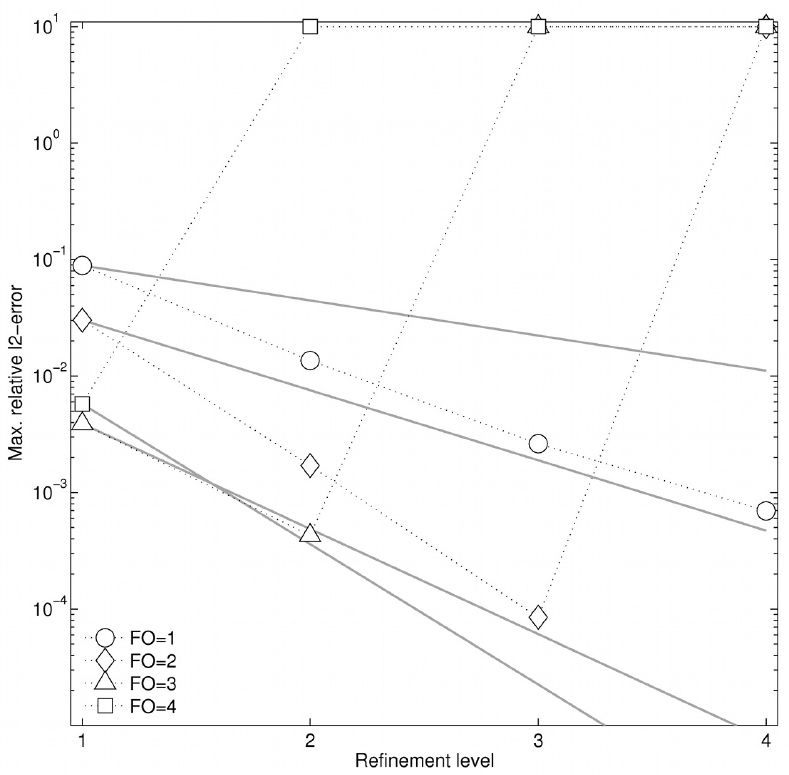}
\caption{Spatial order of convergence for the wave equation integrated until $T=100$. Shown is the maximum $l_2$-error over all generated outputs. The error is capped at 
$10^{1}$, so depicted error values of $10$ correspond to unstable runs.\label{fig:WvEq-SpatConv}}
\end{figure}
\par Figure~\ref{fig:WvEq-SpatConv} shows the maximum $l_2$-error over all generated outputs versus the size of the elements 
of the  employed mesh for finite elements of order one to four. In order to obtain a readable plot, the error is capped at 
$10^{1}$ and values of $10$ correspond to unstable runs. While first order elements are stable on all five meshes, showing 
again better than expected convergence, for higher order elements and fine meshes, the method becomes unstable. Second 
order elements are unstable only on the finest mesh while third and fourth order elements already become unstable after 
three or two refinement steps. Note that on the coarser meshes where the method is stable, the expected or better decay rates of 
the error are observed. Similar behavior is found for the Klein-Gordon equation, but not documented here.
\begin{figure}
	\centering
	\includegraphics[scale=1]{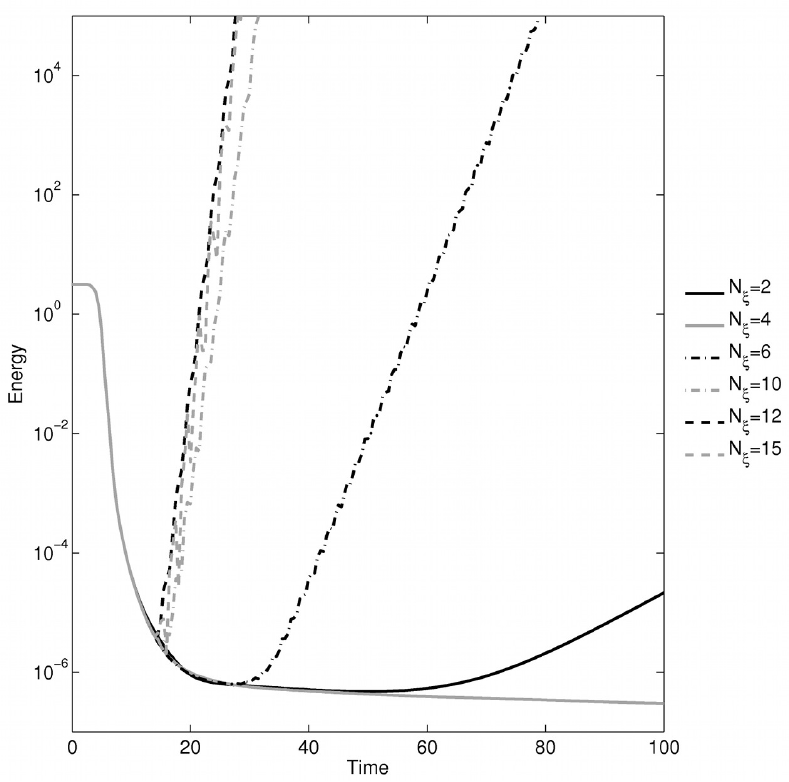}
	\caption{Energy over time for the wave equation for different values of $\Nxi$. The refinement level of the mesh is two and the order of the used finite elements is 
four.\label{fig:WvEq-ErrorVsTime}}
\end{figure}
\par Figure \ref{fig:WvEq-ErrorVsTime} shows the energy of the discrete solution over time for different values of $\Nxi$ on the finest mesh for finite 
elements of order two. The simulations are stable for $\Nxi=4$. They are also stable until about $t=20$ for
 $\Nxi=6, 10, 12, 15$, but an exponential instability occurs after this point in time. The instability  also occurs when the simulation is run on a circular domain. Note that the employed integration scheme is A-stable, so that the instability on a finer mesh is not arising from a violation of some CFL-type stability limit.

\section{Conclusions}
The pole condition approach to transparent boundary conditions, derived in~\cite{RuprechtEtAl2008} for the time-dependent, one-dimensional case, is extended to time-dependent two-dimensional problems. The pole condition identifies in- and outgoing modes 
by associating them with poles of the spatial Laplace transform in the complex plane. The complex plane is then divided into two half-planes, 
$\mathbb{C}_{\rm in}$ and $\mathbb{C}_{\rm out}$, containing the poles corresponding to incoming and outgoing modes respectively. To suppress modes traveling from the
exterior into the computational domain, the Laplace transform is required to be analytic in $\mathbb{C}_{\rm in}$. In
order to obtain a numerically implementable formulation, $\mathbb{C}_{\rm in}$ is mapped to the unit circle by a conformal
M\"obius transformation. The Laplace transform is then extended in a power series on the unit circle with the 
coefficients of the expansion being connected to the interior degrees of freedom on the boundary. Truncating the series 
after a finite number of terms yields an approximate and implementable TBC.
\par
Numerical examples are presented in order to investigate the performance of the pole condition approach:
As in the 1D-case, the considered generic PDE contains different well-known equations for specific choices of parameters.
Excellent results are obtained for Schr\"odinger's equation and the drift-diffusion equation: The presented numerical  experiments demonstrate that the convergence order of finite elements up to order six is retained and also that the convergence order of the temporal discretization is not affected if sufficiently many coefficients are
used for the boundary condition. Further it is shown that the error introduced by
the approximate boundary condition decays super-algebraically as the number of coefficients in the expansion of the Laplace
transform increases. For the drift-diffusion equation, a small minimal number of coefficients was found to be required to reach the regime of super-algebraic error decay.
\par
Unfortunately, in contrast to the one-dimensional case, the approach exhibits instabilities for the two-dimensional wave and Klein-Gordon equation if using finite elements of order two or higher. Hence 
in the present form the pole condition is of limited use for these second order hyperbolic equations. A further investigation of the instability and hopefully a remedy will be subject of future research.

\bibliographystyle{siam}
\bibliography{Mathe,NanoOptics,PoleCondition}

\end{document}